\documentclass[12pt]{amsart}

\usepackage{amsmath,amssymb,amscd,amsthm,amsxtra}
\usepackage[dvips]{graphics,epsfig}

\newtheorem{theorem}{Theorem}

\newtheorem{definition}{Definition}
\newtheorem{corollary}{Corollary}

\newcommand{\q}{\quad}

\newcommand{\nf}{\infty}

\newcommand{\al}{\alpha}
\newcommand{\be}{\beta}

\newcommand{\ve}{\varepsilon}

\newcommand{\rn}{{\mathbf R}^n}

\newcommand{\lab}{\label}

\newcommand{\f}{\frac}

\begin{document}

\title
{Multilinear interpolation between adjoint operators}

\author{Loukas Grafakos}

\address{
Loukas Grafakos\\
Department of Mathematics\\
University of Missouri\\
Columbia, MO 65211, USA}

\email{loukas@math.missouri.edu}

\author{Terence Tao}

\address{
Terence Tao\\
Department of Mathematics\\
University of California, Los Angeles\\
Los Angeles, CA 90024, USA}

\email{tao@math.ucla.edu}

\thanks{Grafakos is supported by the NSF.  Tao is a Clay Prize Fellow and is 
supported by a grant from the Packard Foundation.}

\date{\today}



\begin{abstract}
Multilinear interpolation is a powerful tool used in obtaining strong 
type boundedness for a variety of operators assuming only a finite set 
of restricted weak-type estimates.  A typical situation occurs when one 
knows that a multilinear operator satisfies a weak $L^q$ estimate for a 
single index $q$ (which may be less than one) and that all the adjoints of the 
multilinear operator are of similar nature, and thus they also
 satisfy the same weak $L^q$ estimate.  Under this assumption,  in this expository note 
we give a general multilinear interpolation theorem  
which   allows one to obtain 
strong type boundedness for the operator (and all of its adjoints)
for a large set of exponents.  
The key point in the applications we discuss is that the interpolation theorem can 
handle the case $q \leq 1$.  When $q > 1$, weak $L^q$ has a predual, and
such strong type boundedness 
can be easily obtained by duality and multilinear 
interpolation (c.f. \cite{BL}, \cite{GK}, \cite{janson}, \cite{MTT}, 
\cite{strichartz}). 
\end{abstract}

\maketitle

\section{Multilinear operators}\label{multi-sec}

We begin by setting up some notation for multilinear operators. Let $m 
\geq 1$ be an integer.  We suppose that for 
$0\le j\le m$, $(X_j,\mu_j)$ are measure spaces endowed with positive 
measures $\mu_j$.  We assume that $T$ is an $m$-linear operator of the 
form
$$ T(f_1, \ldots, f_m)(x_0) := \int\ldots\int K(x_0, \ldots, x_m) 
\prod_{i=1}^m f_i(x_i) \ d\mu_i(x_i)$$
where $K$ is a complex-valued locally integrable function on $X_0 
\times \ldots \times X_m$ and $f_j$ are simple functions on $X_j$.  We shall 
make the technical assumption that $K$ is bounded and is supported on a 
product set $Y_0 \times \ldots \times Y_m$ where each $Y_j \subseteq 
X_j$ has finite measure.  Of course, most interesting operators (e.g. 
multilinear singular integral operators) do not obey this condition, but 
in practice one can truncate and/or mollify the kernel of a singular 
integral to obey this condition, apply the multilinear interpolation 
theorem to the truncated operator, and use a standard limiting argument to 
recover estimates for the untruncated operator.

One can rewrite $T$ more symmetrically as an $m+1$-linear form 
$\Lambda$ defined by
$$ \Lambda(f_0, f_1, \ldots, f_m) := \int\ldots\int K(x_0, \ldots, x_m) 
\prod_{i=0}^m f_i(x_i)\ d\mu_i(x_i).$$
One can then define the $m$ adjoints $T^{*j}$ of $T$ for $0 \leq j \leq 
m$
by duality as
$$ \int f_j(x_j) T^{*j}(f_1, \ldots, f_{j-1}, f_0, f_{j+1}, \ldots, 
f_m)(x_j)\ d\mu_j(x_j) := \Lambda(f_0, f_1, \ldots, f_m).$$
Observe that $T = T^{*0}$.

We are interested in the mapping properties of $T$ from the product of spaces
$L^{p_1}(X_1,\mu_1) \times \ldots \times L^{p_m}(X_m, \mu_m)$ into $L^{p_0}(X_0, \mu_0)$ 
for various exponents $p_j$, and more generally for the adjoints 
$T^{*j}$ of $T$.  Actually, it will be more convenient to work with the 
$(m\!+\! 1)$-linear form $\Lambda$, and with the  tuple of   reciprocals $(1/p'_0, 1/p_1, 
\ldots, 1/p_m)$ instead of the exponents $p_j$ directly.  (Here we adopt the 
usual convention that $p'$ is defined by $1/p' + 1/p := 1$ even when $0<p<1$; this notation 
is taken from Hardy, Littlewood and P\'olya.)  

Recall the definition of the weak Lebesgue space $L^{p,\infty}(X_i, \mu_i)$ for $0 < p < 
\infty$ by
$$ \| f \|_{L^{p,\infty}(X_i, \mu_i)} := \sup_{\lambda > 0} \lambda 
\mu_i(\{x_i \in X_i: |f(x_i)| \geq \lambda \})^{1/p}.$$
We also define $L^{\infty,\infty} = L^\infty$.
If $1 < p < \infty$, we define the restricted Lebesgue space 
$L^{p,1}(X_i, \mu_i)$ by duality as
$$ \| f \|_{L^{p,1}(X_i, \mu_i)} := \sup \{ \big|\int f(x_i) g(x_i)\ 
d\mu_i(x_i) \big|: g \in L^{p,\infty}(X_i,\mu_i), \| g\|_{L^{p,\infty}(X_i,\mu_i)} 
\leq 1 \}.$$
We also define $L^{1,1} = L^1$.  This definition is equivalent to the 
other standard definitions of $L^{p,1}(X_i, \mu_i)$ up to a constant 
depending on $p$.

\begin{definition}\label{tuple-def}
Define a \emph{tuple} to be a collection of $m+1$ numbers $\alpha = 
(\alpha_0, \ldots, \alpha_m)$ such that $-\infty < \alpha_i \leq 1$ for 
all $0 \leq i \leq m$, such that such that $\alpha_0 + \ldots + \alpha_m 
= 1$, and such that at most one of the $\alpha_i$ is non-positive. 
If for all $j\in \{0,1,2,\ldots , m\}$  we have $0 
< \alpha_j < 1$, we say that the tuple $\alpha$ is \emph{good}. Otherwise 
there is exactly one $a_i$ such that $a_i\le 0$ and we say 
that the tuple $\alpha$ is \emph{bad}. The smallest number $j_0$ for which 
the $\min\limits_{0\le j\le m}\al_j$ is attained for a tuple 
$\al$ is called the \emph{bad index of the tuple}.
\end{definition}

If $\alpha$ is a good tuple and $B > 0$, we say that $\Lambda$ is 
\emph{of strong-type $\alpha$ with bound $B$} if we have the multilinear 
form estimate
$$ |\Lambda(f_0, \ldots, f_m)| \leq B \prod_{i=0}^m \| f_i 
\|_{L^{1/\alpha_i}(X_i, \mu_i)}$$
for all simple functions $f_0, \ldots, f_m$.  By duality, this is 
equivalent to the multilinear operator estimate
$$ \| T(f_1, \ldots, f_m) \|_{L^{1/(1-\alpha_0)}(X_0, \mu_0)}
\leq B \prod_{i=1}^m \|f_i \|_{ L^{1/\alpha_i}(X_i, \mu_i)}$$
or more generally
$$ \| T^{*j}(f_1, f_{j-1},f_0,f_{j+1},\ldots, f_m) 
\|_{L^{1/(1-\alpha_j)}(X_j, \mu_j)}
\leq B \prod_{\substack{0 \leq i \leq m\\ i \neq j}} \|f_i \|_{ L^{1/\alpha_i}(X_i, 
\mu_i)}$$
for $0 \leq j \leq m$.

If $\alpha$ is a tuple with bad index $j$, we say that $\Lambda$ is 
\emph{of restricted weak-type $\alpha$ with bound $B$} if we have the 
estimate
$$ \| T^{*j}(f_1, f_{j-1},f_0,f_{j+1},\ldots, f_m) 
\|_{L^{1/(1-\alpha_j),\infty}(X_j, \mu_j)}
\leq B \prod_{\substack{0 \leq i \leq m\\ i \neq j}} \|f_i \|_{ 
L^{1/\alpha_i,1}(X_i, \mu_i)}$$
for all simple functions $f_i$.    In view of  duality, if $\al$ is a good index, 
then the choice of the index $j$ above is irrelevant. 

\section{The interpolation theorem}\label{interpolation-sec}

We have the following interpolation theorem for restricted weak-type 
estimates, inspired by  \cite{MTT}:

\begin{theorem}\label{rwt-interp}  Let $\alpha^{(1)}, \ldots, 
\alpha^{(N)}$ be tuples for some $N > 1$, and let $\alpha$ be a good tuple such 
that $\alpha = \theta_1 \alpha^{(1)} + \ldots + \theta_N \alpha^{(N)}$, 
where $0 \leq \theta_s \leq 1$ for all $1 \leq s \leq N$ and $\theta_1 
+ \ldots + \theta_N = 1$.

Suppose that $\Lambda$ is of restricted weak-type $\alpha^{(s)}$ with 
bound $B_s > 0$ for all $1 \leq s \leq N$.  Then $\Lambda$ is of 
restricted weak-type $\alpha$ with bound $C \prod_{s=1}^N B_s^{\theta_s}$, 
where $C > 0$ is a constant depending on $ \alpha^{(1)}, \ldots, 
\alpha^{(N)}, \theta_1, \ldots, \theta_N$.
\end{theorem}

\begin{proof}
Since $\alpha$ is a good tuple, it suffices by duality to prove the 
multilinear form estimate 
$$ | \Lambda(f_0, \ldots, f_m) | \leq C (\prod_{s=1}^N B_s^{\theta_s}) 
\prod_{i=0}^m \| f_i \|_{L^{1/\alpha_i, 1}(X_i, \mu_i)}.$$
We will let the constant $C$ vary from line to line.  
For $1 < p < \infty$, the unit ball of $L^{p,1}(X_i, \mu_i)$ is 
contained in a constant multiple of the convex hull of the normalized 
characteristic functions $\mu_i(E)^{1/p} \chi_E$ (see e.g. 
\cite{sadosky:interp}) it suffices to prove the above estimate for characteristic 
functions:
$$ | \Lambda(\chi_{E_0}, \ldots, \chi_{E_m}) | \leq C (\prod_{s=1}^N 
B_s^{\theta_s}) \prod_{i=0}^m \mu_i(E_i)^{\alpha_i}.$$
We may of course assume that all the $E_i$ have positive finite 
measure.
Let $A$ be the best constant such that
\begin{equation}\label{a-bound}
| \Lambda(\chi_{E_0}, \ldots, \chi_{E_m}) | \leq A (\prod_{s=1}^N 
B_s^{\theta_s}) \prod_{i=0}^m \mu_i(E_i)^{\alpha_i}
\end{equation}
for all such $E_j$; by our technical assumption on the kernel $K$ we 
see that $A$ is finite.  Our task is to show that $A \leq C$.

Let $\ve > 0$ be chosen later.  We may find $E_0, \ldots, E_m$ of 
positive finite measure such that
\begin{equation}\label{a-eps}
| \Lambda(\chi_{E_0}, \ldots, \chi_{E_m}) | \geq (A-\ve)Q,
\end{equation}
where we use $0 < Q < \infty$ to denote the quantity
$$ Q := (\prod_{s=1}^N B_s^{\theta_s}) \prod_{i=0}^m 
\mu_i(E_i)^{\alpha_i}
= \prod_{s=1}^N (B_s \prod_{i=0}^m 
\mu_i(E_i)^{\alpha^{(s)}_i})^{\theta_s}.$$
Fix $E_0, \ldots, E_m$.  From the definition of $Q$
we see that there exists $1 \leq s_0 \leq N$ such that
\begin{equation}\label{q-bound}
B_{s_0} \prod_{i=0}^m \mu_i(E_i)^{\alpha^{(s_0)}_i} \leq Q.
\end{equation}
Fix this $s_0$, and let $j$ be the bad index of $\alpha^{(s_0)}$.  Let 
$F$ be the function
$$ F := T^{*j}(\chi_{E_1}, \ldots, \chi_{E_{j-1}}, \chi_{E_0}, 
\chi_{E_{j+1}}, \ldots, \chi_{E_m}).$$
Since $\Lambda$ is of restricted weak-type $\alpha^{(s_0)}$ with bound 
$B_{s_0}$, we have from \eqref{q-bound} that
\begin{equation}\label{12345}
 \| F \|_{L^{1/(1-\alpha^{(s_0)}_j), \infty}(X_j, \mu_j)} \leq 
B_{s_0} \prod_{\substack{0 \leq i \leq m\\ i \neq j}} \mu_i(E_i)^{\alpha^{(s_0)}_i}
 \leq Q \mu_j(E_j)^{-\alpha^{(s_0)}_j}.
\end{equation}
In particular if we define the set
\begin{equation}\label{epj-def}
E'_j := \{ x_j \in E_j: |F(x_j)| \geq 2^{1-\alpha^{(s_0)}_j} Q 
\mu_j(E_j)^{-1} \}  
\end{equation}
then (\ref{12345}) implies that 
\begin{equation}\label{half}
\mu_j(E'_j) \leq \tfrac{1}{2} \mu_j(E_j).
\end{equation}
By construction of $E'_j$ we have
$|\int \chi_{E_j \backslash E'_j}(x_j) F(x_j)\ d\mu_j(x_j)| \leq 2^{1-\al_j^{(s_0)}} Q$, 
or equivalently that
$$ |\Lambda(\chi_{E_0}, \ldots, \chi_{E_{j-1}}, \chi_{E_j \backslash 
E'_j}, \chi_{E_{j+1}}, \ldots, \chi_{E_m})| \leq CQ.$$
On the other hand, from \eqref{a-bound} and \eqref{half} we have
$$ |\Lambda(\chi_{E_0}, \ldots, \chi_{E_{j-1}}, \chi_{E'_j}, 
\chi_{E_{j+1}}, \ldots, \chi_{E_m})| \leq 2^{-\alpha_j} A Q.$$
Adding the two estimates and using \eqref{a-eps} we obtain
$CQ + 2^{-\alpha_j} A Q \leq (A-\ve) Q$.
Since $\alpha$ is good, we have $\alpha_j > 0$.  The claim $A < C$ then 
follows by choosing $\ve$ sufficiently small.
\end{proof}

>From the multilinear Marcinkiewicz interpolation theorem (see e.g. 
Theorem 4.6 of \cite{GK}) we can obtain strong-type estimates at a good 
tuple $\alpha$ if we know restricted weak-type estimates for all tuples in 
a neighborhood of $\alpha$.  From this and the previous theorem we obtain

\begin{corollary}\label{interp-cor} Let $\alpha^{(1)}, \ldots, 
\alpha^{(N)}$ be tuples for some $N > 1$, and let $\alpha$ be a good tuple in 
the interior of the convex hull of $\alpha^{(1)}, \ldots, \alpha^{(N)}$.
Suppose that $\Lambda$ is of restricted weak-type $\alpha^{(s)}$ with 
bound $B > 0$ for all $1 \leq s \leq N$.  Then $\Lambda$ is of 
strong-type $\alpha$ with bound $C B$, where $C > 0$ is a constant 
depending on $\alpha, \alpha^{(1)}, \ldots, 
\alpha^{(N)}$.
\end{corollary}

By interpolating this result with the restricted weak-type estimates on 
the individual $T^{*j}$, one can obtain some strong-type estimates for 
$T^{*j}$ mapping onto spaces $L^p(X_j, \mu_j)$ where $p$ is possibly 
less than or equal to 1.  By duality one can thus get some estimates 
where some of the functions are in $L^\infty$.  However it is still an open 
question whether one can get the entire interior of the convex hull of 
$\alpha^{(1)}, \ldots, \alpha^{(N)}$ this way\footnote{In \cite{MTT} 
this was achieved, but only after strengthening the hypothesis of 
restricted weak-type to that of ``positive type''.  Essentially, this requires 
the set $E'_j$ defined in \eqref{epj-def} to be stable if one replaces 
the characteristic functions $\chi_{E_i}$ with arbitrary bounded 
functions on $E_i$.}.
 
\section{Applications}

We now pass to three applications. The first application is to re-prove 
an old result of Wolff \cite{wolff}: if $T$ is a linear operator such 
that $T$ and its adjoint $T^*$ both map $L^1$ to $L^{1,\infty}$, then 
$T$ is bounded on $L^p$ for all $1 < p < \infty$ (assuming that $T$ can 
be approximated by truncated kernels as mentioned in the introduction).  
Indeed, in this case $\Lambda$ is of restricted weak-type $(1,0)$ and 
$(0,1)$, and hence of strong-type $(1/p,1/p')$ for all $1 < p < \infty$ 
by Corollary \ref{interp-cor}.

The next application involves the multilinear 
Calder\'on-Zygmund singular integral operators on $\rn\times\dots
\times \rn = (\rn)^m$ defined by
$$
T(f_1, \ldots , f_m)(x_0) :=\lim_{\ve\to 0}  \idotsint\limits_{ 
\sum\limits_{j,k}|x_k-x_j
| \ge \ve} K(x_0,x_1,\ldots , x_m)f_1(x_1)\ldots f_m(x_m)\, 
\,dx_1\ldots dx_m , 
$$
where $ |K(\vec x)|\le C (\sum_{j,k=0}^m|x_k-x_j|)^{-nm}$, 
$ |\nabla K(\vec x)|\le C (\sum_{j,k=0}^m|x_k-x_j|)^{-nm-1}$, and 
$\vec x= (x_0,x_1, \dots , x_m)$.  
These integrals have been extensively studied by Coifman and Meyer 
\cite{CM1},\cite{CM2},\cite{CM3} and recently by Grafakos and Torres 
\cite{GT}. 
It was shown in \cite{GT} and also by Kenig and Stein \cite{KS} 
(who considered the case $n\! =\! 1$,  $m\! =\! 2$) that if such 
operators map 
$L^{q_1} \times \dots \times L^{q_m} $ into $L^{q,\nf} $ for only one 
$m$-tuple of indices, then they must map $L^1\times\dots \times L^1$ 
into $L^{1/m,\nf}$.  
Since the adjoints of these operators satisfy similar boundedness 
properties, we see that the corresponding form $\Lambda$ is of 
restricted weak-type $(1-m, 1,\ldots,1)$, and similarly for permutations.  It 
then follows\footnote{Strictly speaking, we have to first fix $\ve$, and 
truncate the kernel $K$ to a compact set, before applying the Theorem, 
and then take limits at the end.  We leave the details of this standard 
argument to the reader.  A similar approximation technique can be 
applied for the bilinear Hilbert transform below.}  from Corollary 
\ref{interp-cor} that $T$ maps  
$L^{p_1} \times \dots \times L^{p_m} $ into $L^{p} $ for all $m$-tuples 
of indices
with\footnote{The convex hull of the permutations of $(1-m,1,\ldots,1)$ 
is the tetrahedron of points $(x_0, \ldots, x_m)$ with $x_0 + \ldots + 
x_m = 1$ and $x_i \leq 1$ for all $0 \leq i \leq m$, so in particular 
the points $(1/p_1, \ldots, 1/p_m)$ described above fall into this 
category.} $1<p_j<\nf$ with $\f{1}{p_1}+\dots + \f{1}{p_m}=\f{1}{p}$ and $p 
> 1$.  The condition $p>1$ can be removed by further interpolation with 
the $L^1 \times \ldots \times L^1 \to L^{1/m}$ estimate. This argument 
simplifies the interpolation proof used in \cite{GT}.  

Our third application involves the bilinear Hilbert transform 
$H_{\al,\be}$ defined by
\begin{equation}\lab{bht} 
H_{\al,\be}(f,g)(x) = \lim_{\ve\to 0} \int_{|t|\ge \ve} 
f(x-\al t)g(x-\be t)\, \f{dt}{ t} \, ,\q  x\in \mathbf R\, .  
\end{equation} 
The proof of boundedness of $H_{\al,\be}$ from $L^2\times L^2$ into
$L^{1,\nf}$   (for example see   \cite{lacey})
 is technically simpler than that of $L^{p_1}\times L^{p_2}$
into $L^p$ when $2<p_1,p_2,p'<\nf$   given in 
Lacey and Thiele \cite{lacey-thiele1}. Since the adjoints of the 
operators $H_{\al,\be}$ are $H_{\al,\be}^{*1} =H_{-\al,\be-\al}$  and 
$H_{\al,\be}^{*2} =H_{ \al-\be,-\be}$ which  are ``essentially'' the 
same 
operators, we can use the single estimate $L^2\times L^2\to L^{1,\nf}$ 
for 
all of these operators to obtain the results in \cite{lacey-thiele1}, 
since the corresponding form $\Lambda$ is then of restricted weak-type 
$(0,1/2,1/2)$, $(1/2,0,1/2)$, and $(1/2,1/2,0)$.  (See also the similar 
argument in \cite{MTT}).  

The operator in \eqref{bht} is in fact bounded in the larger range $1 < 
p_1, p_2 < \infty$, $p > 2/3$ and similarly for adjoints, see 
\cite{lacey-thiele2}.  The interpolation theorem given here allows for a slight 
simplification in the arguments in that paper (cf. \cite{MTT}), 
although one cannot deduce all these estimates solely from the $L^2 \times L^2 
\to L^{1,\infty}$ estimate.




\begin{thebibliography}{43}

\bibitem{BL} J. Bergh and J. L\"ofstr\"om, 
\emph{Interpolation spaces, An introduction}, 
Springer-Verlag, New York, NY 1976. 

\bibitem{CM1} R. R.  Coifman and Y.  Meyer, {\it On commutators of 
singular
integrals and bilinear singular integrals}, Trans.  Amer.  Math.  Soc. 
\textbf{212} (1975), 315--331.

\bibitem{CM2} R. R.  Coifman and Y. Meyer, {\it Commutateurs d' 
int\'egrales
singuli\`eres et op\'erateurs multilin\'eaires}, Ann. Inst. Fourier, 
Grenoble
\textbf{28} (1978), 177--202.

\bibitem{CM3} R. R.  Coifman and Y.  Meyer, {\it Au-del\`a des 
op\'erateurs
pseudo-diff\'erentiels}, Asterisk \textbf{57}, 1978.

\bibitem{GK} L. Grafakos and N. Kalton,
{\it Some remarks on multilinear maps and interpolation},
Math. Ann. {\bf 319} (2001), 151--180. 
 
\bibitem{GT} L. Grafakos and R. Torres,
{\it Multilinear Calder\'on-Zygmund theory},
Advaces in Math.,  to  appear. 


\bibitem{janson} S.  Janson, {\it On interpolation of multilinear 
operators},
in  Function spaces and applications (Lund, 1986),
Lecture Notes in Math. {\bf 1302},
Springer-Verlag, Berlin-New York, 1988.

\bibitem{KS} C.   Kenig and E.  M.  Stein,
{\it Multilinear estimates and fractional integration},
Math. Res. Lett. {\bf 6} (1999), 1--15.

\bibitem{lacey} M.  T.   Lacey,
{\it On the bilinear Hilbert transform},
Doc. Math. {\bf 1998}, Extra Vol. II, 647--656. 

\bibitem{lacey-thiele1} M.  T.   Lacey and C.  M.  Thiele,
{\it $L^p$ bounds for the bilinear Hilbert transform, $2<p<\infty$},
Ann. Math. {\bf 146} (1997), 693--724.

\bibitem{lacey-thiele2} M.  T. Lacey and C.  M.  Thiele,
{\it On Calder\'on's conjecture},  Ann.   Math. {\bf 149} (1999), 
475--496.

\bibitem{MTT} C. Muscalu, C. Thiele, and T. Tao, 
{\it  Multi-linear operators given by singular multipliers}, 
J. Amer. Math. Soc., to appear.

\bibitem{sadosky:interp}
C. Sadosky, \emph{Interpolation of Operators and Singular Integrals},
Marcel Dekker Inc., 1976.
 
\bibitem{strichartz} R.\ Strichartz, {\it A multilinear version of the
Marcinkiewicz interpolation theorem}, Proc.\ Amer.\ Math.\ Soc. {\bf 
21}
(1969), 441--444.

\bibitem{wolff} T. H. Wolff. \emph{A note on interpolation spaces}. 
Harmonic analysis (Minneapolis, Minn., 1981), pp. 199--204, Lecture Notes 
in Math., 908, Springer, Berlin-New York, 1982.

\end{thebibliography}
\end{document}